\documentclass{article}

\usepackage{graphicx} 

\usepackage{latexsym, amscd, amsfonts, eucal, mathrsfs, amsmath, amssymb, amsthm, xypic,xr, makeidx, stmaryrd}

\usepackage{tikz-cd}
\usepackage{bookmark}
\usepackage{makeidx}
\usepackage{geometry}
\usepackage{indentfirst}
\usepackage{enumitem}
\usepackage{hyperref}
\usepackage{mathrsfs}
\usepackage{enumitem}
\usepackage{verbatim}
\usepackage{upgreek}
\usepackage{textcomp}
\usepackage{amsmath}
\usepackage{dsfont}

\newtheorem{thm}{Theorem}[subsection]
\newtheorem{defn}[thm]{Definition}

\newtheorem{prop}[thm]{Proposition}

\newtheorem{example}[thm]{Example}
\newtheorem{remark}[thm]{Remark}

\AfterEndEnvironment{thm}{\noindent\ignorespaces}
\AfterEndEnvironment{defn}{\noindent\ignorespaces}
\AfterEndEnvironment{notation}{\noindent\ignorespaces}
\AfterEndEnvironment{prop}{\noindent\ignorespaces}
\AfterEndEnvironment{lem}{\noindent\ignorespaces}
\AfterEndEnvironment{cor}{\noindent\ignorespaces}
\AfterEndEnvironment{example}{\noindent\ignorespaces}
\AfterEndEnvironment{remark}{\noindent\ignorespaces}

\DeclareMathOperator{\sym}{\mathbb{S}}
\DeclareMathOperator{\X}{{\mathcal{X}}}
\DeclareMathOperator{\OX}{\mathcal{O}_{\X}}
\DeclareMathOperator{\G}{{\mathcal{G}_{\sym} ^{alg}}}
\DeclareMathOperator{\T}{{\mathcal{T}_{\sym} ^{alg}}}
\DeclareMathOperator{\shv}{\mathrm{Shv}}
\DeclareMathOperator{\psh}{\mathrm{PSh}}
\DeclareMathOperator{\spec}{\mathrm{Spec}}
\DeclareMathOperator{\fun}{\mathrm{Fun}}
\DeclareMathOperator{\map}{\mathrm{Map}}

\DeclareMathOperator{\perf}{\mathrm{Perf}}

\DeclareMathOperator{\Sp}{\mathrm{Sp}}
\DeclareMathOperator{\Symp}{\mathrm{Symp}}
\DeclareMathOperator{\symp}{\mathrm{symp}}
\DeclareMathOperator{\idem}{\mathrm{Idem}}
\DeclareMathOperator{\mot}{\mathrm{Mot}}

\title{Some memos on Stable Symplectic Structured Space II: Symplectic motives}
\author{Eita Haibara}
\date{April 2025}

\begin{document}

\maketitle
\begin{abstract}
These memos include the research on $\G$-scheme theory, the definition of symplectic motives over $\G$-schemes and symplectic motivic cohomology. This presents a new research direction.
\end{abstract}
\tableofcontents

\section{Scheme theory}
We assume that the reader reads \cite{haibara1} and we assume that the reader is familiar with higher algebra \cite{ha}. We follow the notation of \cite{haibara1}. This paper is an extension of paper \cite{haibara1}.
\subsection{Scheme theory}

\begin{defn}
An open immersion is a morphism $(j,\beta) : (\mathcal{U},\mathcal{O}_{\mathcal{U}}) \to (\X,\OX)$ of $\G$-schemes such that $j_*: \mathcal{U} \to \X$ is an open embedding of $\infty$-topoi and $\beta: j^*\OX \to \mathcal{O}_{\mathcal{U}}$ is an isomorphism.
\end{defn}
\begin{defn}
    Embedding $\T \hookrightarrow \G$ sends $(M,M^{-\tau})$ to $(M,A)$, where $A$ is constructed from $M^{-\tau}$. Specifically, we define the sheaf $A$ on $M$ by $A(U)= U^{-\tau}$, for each open subset $U \subseteq M$, with restriction maps $\pi : M^{-\tau} \to U^{-\tau}$ provided by the Thom spectrum structure over open embeddings.
\end{defn}

Following example is important in this paper.
\begin{example}
Consider the symplectic manifold $M=\mathbb{CP}^1$, the complex projective line, equipped with the Fubini-Study symplectic form $\omega$. This is a compact, algebraic variety with a natural symplectic structure, making it suitable for our purposes. Define $A = \mathcal{O}_{\mathbb{CP}^1}$ , the structure sheaf of holomorphic functions on $\mathbb{CP}^1$. We can view $\mathcal{O}_{\mathbb{CP}^1}$ as an discrete $E_{\infty}$-ring. Consider $\infty$-topos $\X=\shv(\mathbb{CP}^1)$, $\OX: \G \to \X$ and $\OX$ assigns to $\mathbb{CP}^1$ the sheaf $\mathcal{O}_{\mathbb{CP}^1}$. Thus, the $\G$-scheme is $(\X,\OX)$, representing $(\mathbb{CP}^1,\mathcal{O}_{\mathbb{CP}^1})$ with its symplectic structure. Then consider $\perf(\X)$, these are nice sheaves (like vector bundles and their derived analogs) on $\mathbb{CP}^1$. Since $\mathbb{CP}^1$ is a smooth projective variety, $\perf(\X)$ is equivalent to the derived category of perfect complexes on $\mathbb{CP}^1$. It’s generated by line bundles $\mathcal{O}(n)$ (where $n \in \mathbb{Z}$), like $\mathcal{O}(0)$ (the trivial bundle) and $\mathcal{O}(1)$ (the hyperplane bundle). Define $\Symp_0(\X)$ as the stable $\infty$-category of perfect complexes $E \in \perf(\X)$ equipped with a symplectic form: A non-degenerate pairing $\omega_E : E \otimes E \to \OX$, where the shift $d = 0$ indicates an unshifted symplectic structure(analogous to a classical symplectic vector bundle). Take $E = \mathcal{O}_{\mathbb{CP}^1} \bigoplus \mathcal{O}_{\mathbb{CP}^1}$, a trivial rank-2 vector bundle. Define $\omega_E: E \otimes E \to  \mathcal{O}_{\mathbb{CP}^1}$ by $(f_1,f_2) \otimes (g_1,g_2) \mapsto f_1g_2 - f_2g_1$. This is non-degenerate so $(E,\omega_E)$ is an object in $\Symp_0(\X)$.
\end{example}

Note that: Every vector bundle E on $\mathbb{CP}^1$ splits as a direct sum of line bundles, $E \cong \bigoplus_{i=1}^{r}\mathcal{O}(k_i)$, where $\mathcal{O}(k_i)$ are line bundle with Chern class $k_i \in \mathbb{Z}$, and r is the rank of E. \cite{grothendeick57} For a vector bundle E to admit symplectic form, the rank must be even and the symplectic form $\omega_E : E \otimes E \to \mathcal{O}_{\mathbb{CP}^1}$ identifies E with its dual $E^{\vee}$ in a skew-symmetric way, i.e. $E\cong E^{\vee}$ up to a twist.

\subsection{Morphisms of schemes}
Consider a morphism $f$ that is a closed immersion. Then for each $(M,A) \in \G$, we have [\cite{haibara1}, Definition 2.3.8]:
\newline
(1) $\mathcal{O}_{\mathcal{Y}} \in \mathcal{Y}$, an object in the $\infty$-topos $\mathcal{Y}$.
\newline
(2) $f^{*}\mathcal{O}_{\mathcal{Y}}(M,A) \in \X$, its pullback to $\X$
\newline
(3) $\OX(M,A) \in \X$, the structure sheaf of $\X$ evaluated at $(M,A)$
\newline
(4) $\beta_{(M,A)} : f^{*}\mathcal{O}_{\mathcal{Y}}(M,A) \to \OX(M,A)$, an effective epimorphism.
\newline
\begin{defn}
The map $\beta_{(M,A)}$ is a morphism in $\X$. Assuming that $\OX(M,A)$ and $f^*\mathcal{O}_{\mathcal{Y}}(M,A)$ are objects in a stable subcategory of $\X$, such as the category of spectrum objects $\Sp(\X)$, we define the kernel of $\beta_{(M,A)}$ as the fiber:
\begin{equation*}
    \mathrm{ker}(\beta_{(M,A)}) = \mathrm{fib}(\beta_{(M,A)}) = \mathrm{fib}(f^{*}\mathcal{O}_{\mathcal{Y}}(M,A) \to \OX(M,A)), 
\end{equation*}
where the fiber is the homotopy pullback in this stable subcategory, computed as the pullback along $f^{*}\mathcal{O}_{\mathcal{Y}}(M,A) \to \OX(M,A) \leftarrow 0$, with $0$ being the zero object in the stable subcategory.
\end{defn}
We interpret the kernel of $f^*$ as the collection $\mathcal{I} = [\mathrm{ker}(\beta_{(M,A)})]_{(M,A) \in \G}$, a family of objects in $\X$ indexed by $\G$.

\begin{prop}
    Since $\beta_{(M,A)}$ is an effective epimorphism, $\OX(M,A)$ is the cofiber of $\mathcal{I}_{(M,A)} \to f^{*}\mathcal{O}_{\mathcal{Y}}(M,A)$. If $\mathcal{O}_{\mathcal{Y}}(M,A)$ and $\OX(M,A)$ are connective objects and $\beta_{(M,A)}$ induce surjection on $\pi_0$, then $\mathcal{I}_{(M,A)}$ is connective.
\end{prop}

\begin{defn}
An object (e.g. a sheaf of spectra or an $E_{\infty}$-ring on the $\infty$-topos $\X$) is said to be nilpotent if, for every open subset $U \subseteq  \X$, the section $\mathcal{F}(U)$ has a finite Postnikov tower. That is, $\mathcal{F}(U)$ is connective (no negative homotopy groups) and bounded (only finitely many positive homotopy groups). 
\end{defn}
This ensures nilpotency is a local property, consistent with the sheaf-theoretic framework.
\newline
\\
Consider a morphism of $\G$-schemes $f: (\X,\OX) \to (\mathcal{Y},\mathcal{O}_{\mathcal{Y}})$, defined by [\cite{haibara1},  Definition 2.3.5]. This induces a functor between the associated categories: $f^{*}:\Symp_d(\mathcal{Y}) \to \Symp_d(\mathcal{X})$. 
\begin{prop}
    Since $f: \X \to \mathcal{Y}$ is a morphism of $\G$-schemes, the induced functor $f^{*}:\Symp_d(\mathcal{Y}) \to \Symp_d(\mathcal{X})$ is exact. 
    \newline
    \newline
    Proof) Define: $f^*(E,\omega_E) = (f^*E,f^*\omega_E)$ remains in $\Symp_d(\mathcal{X})$, where $f^*E$ is the pullback of the complex $E$ to $\X$, $f^*\omega_E$ is the induced symplectic form on $f^*E$, given by the composition: 
    \begin{equation*}
        f^*E \otimes f^*E \xrightarrow[]{f^*\omega_E}f^*(\mathcal{O}_{\mathcal{Y}}[d])\xrightarrow[]{\beta[d]}\OX[d]
    \end{equation*}
    It preserves the zero object, mapping $(0,\omega_0)$ to $(0,0)$. An exact triangle in $\Symp_d(\mathcal{Y})$ is of the form: 
    \begin{equation*}
        A\xrightarrow[]{u}B\xrightarrow[]{v}C\xrightarrow[]{w}\Sigma A
    \end{equation*}
    where $A = (A,\omega_A), B = (B,\omega_B), C = (C,\omega_C)$, and $\Sigma A$ is the suspension of $A$(a shift in the stable structure), with maps $u,v,w$ preserving symplectic forms via homotopies. Apply $f^*$:
    \begin{equation*}
        f^*A\xrightarrow[]{f^*u}f^*B\xrightarrow[]{f^*v}f^*C\xrightarrow[]{f^*w}f^*(\Sigma A)
    \end{equation*}
    Since $f^*: \perf(\mathcal{Y}) \to \perf(\X)$ is exact, it commutes with suspension: $f^*(\Sigma A) \cong \Sigma(f^*A).$ In $\perf(\mathcal{Y})$ this is an exact triangle, and $f^*$ preserves this structure because it is exact as a functor between stable $\infty$-categories. Since $u: A \to B$ preserves the symplectic form (i.e. $(u \otimes u)  \circ  \omega_A \cong \omega_B  \circ u$ up to homotopy), applying $f^*$ to this homotopy uses the functoriality of $f^*$ on tensor products and compositions, preserving the equivalence in $\Symp_d(\X)$. Similarly for $v$ and $w$. Thus, $f^*A \to f^*B \to f^*C \to f^*(\Sigma A)$ is an exact triangle in $\Symp_d(\X)$. Since $\Symp_d(\mathcal{Y})$ and $\Symp_d(\mathcal{X})$ are stable, exact triangles generate all finite limits and colimits. 
\end{prop}

\section{Symplectic K-Theory}
In this section, we assume that the reader reads \cite{Blumberg_2013}.
\subsection{Relative Symplectic K-Theory}

\begin{defn}
Consider the fiber sequence of spectra: 
\begin{equation*}
    K(\Symp_d(f)) \to K(\Symp_d(\mathcal{Y})) \to K(\Symp_d(\X))
\end{equation*}
Here, $K(\Symp_d(f))$ is the relative K-theory spectrum associated with the morphism $f$, analogous to $K(A,I)$ in the classical setting, where $I$ is the ideal defining the quotient $A \to A/I$. 
\end{defn}
For a small idempotent-complete stable $\infty$-category $A$, category of compact spectra $S_{\infty}^{\omega}$, connective K-theory($K$) and non-connective K-theory($\mathbb{K}$), we have [\cite{Blumberg_2013}, Theorem 1.3.]: 
\begin{equation*}
    \map(\mathcal{U}_{add}(S_{\infty}^{\omega}),\mathcal{U}_{add}(A)) \cong K(A)
\end{equation*}
\begin{equation*}
    \map(\mathcal{U}_{add}(S_{\infty}^{\omega}),\mathcal{U}_{loc}(A)) \cong \mathbb{K}(A)
\end{equation*}
These equivalences allow K-theory to be studied via mapping spectra in $\mathcal{M}_{add}$ or $\mathcal{M}_{loc}$.
\newline
\\
Since $\Symp_d(\X)$ is may not split all idempotents due to symplectic constraints, we will use the $\idem(\Symp_d(\X))$.
\begin{prop}
    Assuming $\Symp_d(\X)$ is small(more precisely, that we are working in a universe where $\Symp_d(\X)$ is small relative to a chosen Grothendieck universe) it belongs to $\mathrm{Cat}_{\infty}^{ex}$, and its idempotent completion $\idem(\Symp_d(\mathcal{X}))$ is a small idempotent-complete stable $\infty$-category. The connective K-theory satisfies: 
    \begin{equation*}
        K(\idem(\Symp_d(\mathcal{X}))) \cong \map(\mathcal{U}_{add}(S_{\infty}^{\omega}),\mathcal{U}_{add}(\idem(\Symp_d(\X))))
    \end{equation*}
\end{prop}
\begin{prop}
    On the above assumptions, let $f: \X \to \mathcal{Y}$ be an open immersion. Consider the exact sequence of stable $\infty$-categories:
    \begin{equation*}
        \Symp_d(f) \to \Symp_d(\mathcal{Y}) \xrightarrow[]{f^*} \Symp_d(\X).
    \end{equation*}
    This sequence is exact in $\mathrm{Cat}_{\infty}^{ex}$. \\
    Specifically, $\Symp_d(f)$ is the full subcategory of $\Symp_d(\mathcal{Y})$ consisting of objects that map to zero in $\Symp_d(\X)$ under the functor $f^*$. This ensures that $\Symp_d(f) \to \Symp_d(\mathcal{Y})$ is fully faithful, and $\Symp_d(\X)$ is equivalent to the quotient $\Symp_d(\mathcal{Y})/\Symp_d(f)$. Then non-connective K-theory of spectra $\mathbb{K}(\idem(\Symp_d(f)))$ is the fiber of $\mathbb{K}(\idem(\Symp_d(\mathcal{Y}))) \to \mathbb{K}(\idem(\Symp_d(\mathcal{X})))$. And this fiber corresponds to:
    \begin{equation*}
    \map(\mathcal{U}_{loc}(S_{\infty}^{\omega})\text{, fiber of }\mathcal{U}_{loc}(\idem(\Symp_d(\mathcal{Y}))) \to \mathcal{U}_{loc}(\idem(\Symp_d(\mathcal{X})))).
    \end{equation*}
\end{prop}

\subsection{Applying the Topological Dennis Trace}
\begin{remark}
    For a stable $\infty$-category $\mathcal{A}$, the Grothendieck group $K_0(A)$ is an abelian group defined as follows:
    \newline
    Generators: Isomorphism classes $[A]$ of objects $A \in \mathcal{A}$.
    \newline
    Relations: For every exact triangle $A \to B \to C \to \sum A$, impose the relation $[B] = [A] + [C]$.
\end{remark}
A stable $\infty$-category $A$ in $\mathrm{Cat}_{\infty}^{perf}$(idempotent-complete small stable $\infty$-categories) is dualizable if it is smooth and proper.[\cite{Blumberg_2013}, Theorem 3.7.]  For such $A$, its dual is $A^{op}$, and it's K-theory $K(A)$ admits topological Dennis trace map to $THH(A)$. 

\begin{prop}
    If $\idem(\Symp_d(\X))$ is smooth and proper, then it is dualizable, and its K-theory $K(\idem(\Symp_d(\X)))$ admits a topological Dennis trace map to $THH(\idem(\Symp_d(\X)))$.
\end{prop}
\newline
Let us revisit [Example 1.1.3.] for illustration.
\begin{example}
    Since $\X = \shv(\mathbb{CP}^1)$, we have $\perf(\X) \cong \perf(\mathbb{CP}^1)$, and perfect complexes on $\mathbb{CP}^1$ are equivalent to vector bundles(up to homotopy). Since $\mathbb{CP}^1$ is a smooth projective variety, objects in $\Symp_0(\X)$ can be thought of as symplectic vector bundles on $\mathbb{CP}^1$. For $\omega_e$ to be non-degenerate and skew-symmetric: The pairing $\mathcal{O}(k_i) \otimes \mathcal{O}(k_j) \to \mathcal{O}$ is given by $Hom(\mathcal{O}(k_i),\mathcal{O}(-k_j)) = H^0(\mathbb{CP}^1,\mathcal{O}(-k_j -k_i))$. This is non zero(i.e. $\mathbb{C}$) only when $-k_j -k_i = 0$, so $k_j = -k_i$ (note that: $\mathrm{dim} H^0(\mathbb{CP}^1,\mathcal{O}(0)) = 1.)$. Thus, $E$ must be pair line bundles of opposite degrees, leading to: 
    \begin{equation*}
        E \cong \bigoplus_{I=1}^n(\mathcal{O}(k_i) \bigoplus \mathcal{O}(-k_i)),
    \end{equation*}
    where $k_i \geq 0$ (since $\mathcal{O}(k_i) \bigoplus \mathcal{O}(-k_i) \cong \mathcal{O}(-k_i) \bigoplus \mathcal{O}(k_i)$), and $\omega_E$ pairs each $\mathcal{O}(k_i)$ with $\mathcal{O}(-k_i)$. (Indecomposables: $(\mathcal{O}(k_i) \bigoplus \mathcal{O}(-k_i),\omega_k)$ for $k_i \geq 0$.)
    Now we recall the algebraic K-theory of $\mathbb{CP}^1$: $K_0(\mathbb{CP}^1)$: Generated by $[\mathcal{O}]$ and $[\mathcal{O}(1)]$, with $K_0 \cong \mathbb{Z} \bigoplus \mathbb{Z}$, corresponding to rank and degree (or Chern class $c_1$). Exact triangles are of the form $(A,\omega_A) \to (B,\omega_B) \to (C,\omega_C) \to \sum(A,\omega_A)$, inherited from $\perf(\mathbb{CP}^1)$. On $\mathbb{CP}^1,$ vector bundles are semi-simple, so exact sequences split, and thus: $(B,\omega_B) \cong (A,\omega_A) \bigoplus (C,\omega_C)$. The Grothendieck group $K_0(\Symp_0(\X))$ for $\X = \shv(\mathbb{CP}^1)$ is defined as generators $[E,\omega_E]$ for all objects, relations $[B] = [A] + [C]$ when $B \cong A \bigoplus C$ and basis $[(\mathcal{O}(k) \bigoplus \mathcal{O}(-k),\omega_k)]$ for $k\geq 0$. So the result is 
    \begin{equation*}
        K_0(\Symp_0(\X)) \cong \bigoplus_{k=0}^{\infty} \mathbb{Z} \cdot [(\mathcal{O}(k) \bigoplus \mathcal{O}(-k),\omega_k)] \cong \mathbb{Z}^{\bigoplus \mathbb{N}}.
    \end{equation*}
   Since $\idem(\Symp_0(\X))$ adds formal summands, and stability ensures these are already accounted for, we have: 
    \begin{equation*}
        K_0(\Symp_0(\X)) \cong K_0(\idem(\Symp_0(\X))).
    \end{equation*}
\end{example}

\section{Symplectic motives over $\G$-schemes}
\subsection{Category of motives}
Let's take stable $\infty$-category of spectrum objects in $\X$, $\Sp(\X)$ from [\cite{ha}, Section 1.4]. In the pregeometry $\T$ underlying $\G$, Thom spectra $M^{-\tau}$ are $E_{\infty}$-ring spectra tied to the stable normal bundle of $M$. These encode symplectic invariants, such as the symplectic form. To incorporate the symplectic structure, we define a global Thom spectrum object $\theta_{\X} \in \Sp(\X)$:
\begin{defn}
   Local definition: For an affine $\G$-scheme $\spec^{\G}(A) = (\X_A, \OX_A)$, where $A = (M,A_M)$, the Thom spectrum $M^{-\tau}$ (associated with the stable normal bundle of symplectic manifold $M$) is an $E_{\infty}$-ring spectrum in $\Sp(\X_A)$. Global Construction: For a general $\G$-scheme $\X$, the global Thom spectrum $\theta_{\X} \in \Sp(\X)$ is defined by gluing the local Thom spectra $M^{-\tau}$ on affine opens $\X_A \hookrightarrow \X$. Specifically, for each open immersion $j: \X_A \to \X$, we have $j^*\theta_{\X} \cong M^{-\tau}$. On overlaps $\X_A \times_{\X} \X_B = \X_C$, the restrictions $j_1^*M^{-\tau} \cong U^{-\tau} \cong j_2^*N^{-\tau}$ hold via the restriction maps $\pi: M^{-\tau} \to U^{-\tau}$ and $pi: N^{-\tau} \to U^{-\tau}$, ensuring consistent gluing.
\end{defn}

\begin{defn}
A symplectic motive over a $\G$-scheme $\X$ is an object in $\Sp(\X)$ equipped with a symplectic structure. Formally, objects are a pair $(E,\omega_E)$, where $E \in \Sp(\X)$ is a spectrum object and $\omega_E: E \otimes E \to \theta_{\X}[d]$ is a map in $\Sp(\X)$(for some integer $d$), where $\theta_{\X}[d]$ is the $d$-fold shift of $\theta_{\X}$. The induced map $\upphi_{\omega_E}: E \to \underline{Hom}(E,\theta_{\X}[d])$, defined by $\upphi_{\omega_E}(e)(e') = \omega_E(e \otimes e')$, is an equivalence in $\Sp(\X)$ (non-degeneracy). $\omega_E$ is skew-symmetric, meaning $\omega_E \circ \tau \cong -\omega_E$, where $\tau: E \otimes E \to E \otimes E$ is the symmetry isomorphism in $\Sp(\X)$, and $-\omega_E: E\otimes E\xrightarrow[]{\omega_E}\theta_{\X}[d]\xrightarrow[]{-1}\theta_{\X}[d]$ with $-1$ the negation map in $\theta_{\X}[d]$.
\end{defn}

\begin{defn}
    The category of symplectic motives over $\X$, denoted $\mot_{\symp}(\X)$, is defined as follows:
    \begin{enumerate}[label=(\arabic*)]
        \item Objects: Pairs $(E,\omega_E)$, as described above.
        \item Morphisms: A morphism $f:(E,\omega_E) \to (F,\omega_F)$ is a map $f: E \to F$ in $\Sp(\X)$ that preserves the symplectic structure up to homotopy. Specifically, there exists a homotopy $h:\omega_F \circ ( f \otimes f) \cong \omega_E$, where both sides map  $E \otimes E \to \theta_{\X}[d]$, ensuring f preserves the symplectic structure up to homotopy.
    \end{enumerate}
\end{defn}

\subsection{Symplectic Nisnevich topology}
\begin{defn}
    $C_{\X}$ is defined as follows:
    \begin{enumerate}[label=(\arabic*)]
        \item Objects: Let $\mathcal{Y}$ be a scheme equipped with $\G$-strucuture, meaning it has an underlying symplectic manifold structure via the geometry $\G$. Pair ($V,\mathcal{O}_V$), where V is an étale scheme over $\mathcal{Y}$, meaning the morphism $V \to \mathcal{Y}$ is étale. $\mathcal{O}_V: 
        \G \to Shv(V)$ is a $\G$-structure on $Shv(V)$, making $(Shv(V),\mathcal{O}_V)$ a $\G$-scheme. There exists a morphism of $\G$-schemes $(Shv(V),\mathcal{O}_V) \to (\X,\OX)$, which is induced by the étale map $V \to \mathcal{Y}$ and the relationship $\X =  Shv(\mathcal{Y}).$
        
        \item Morphisms: A morphism ($V,\mathcal{O}_V) \to (W,\mathcal{O}_W$) is an open embedding $V \hookrightarrow W$ of schemes that is compatible with the $\G$-structures, ensuring the induced map on sheaves preserves the algebraic and symplectic data.
    \end{enumerate}
\end{defn}
\begin{defn}
    Symplectic Nisnevich topology is defined as follows:
    \newline
    \newline
    Covering: A covering in this topology is a family of morphisms $[(U_i,\mathcal{O}_{U_i}) \to (U,\mathcal{O}_U)]$, where each morphism is an open embedding of $\G$-schemes, and the induced maps on the $\G$-structures preserve the derived symplectic data. The family is surjective, meaning the union of the images of the $U_i$'s covers $U$ completely. 
    \end{defn}
This is analogous to the Nisnevich topology in algebraic geometry, where coverings involve étale maps that hit every point, but here we emphasize open embeddings and symplectic compatibility. 
\newline
\newline
For $C_{\X}$ to work as a site, it must have fiber products, so we can compute intersections of open sets. And the topology must satisfy standard axioms: coverings are stable under pullbacks, transitive, and non-empty for every object. These properties hold because open embeddings behave well in topoi, and the symplectic Nisnevich topology is designed to mimic well-behaved topologies like Zariski or Nisnevich in classical settings.
\newline
\newline
Begin with the $\infty$-category of presheaves $\psh(C_{\X}) = \fun(C_{\X}^{op},S)$ where $S$ is the $\infty$-category of spaces. Localized $\psh(C_{\X})$ with respect to the symplectic Nisnevich topology to obtain $\shv(C_{\X})$. This localization enforces the descent condition: sheaves must glue consistently over symplectic Nisnevich coverings.
\begin{defn}
    The class $W$ is defined as the set of maps $\upphi:U \to V$, where $U,V \in C_{\X}$, satisfying, $\upphi$ is a symplectomorphism, and there exists a symplectomorphism $\upphi':V\to U$ such that:
    \newline
    \newline 
    The composition $\upphi' \circ \upphi$ is isotopic to the identity map $id_U$ through a smooth family of symplectomorphisms $\upphi_t:U\to U$ (for $t \in [0,1]$), where $\upphi_0=id_U$ and $\upphi_1 = \upphi' \circ \upphi$.
    \newline
    The composition $\upphi \circ \upphi'$ is isotopic to the identity map $id_V$ through a smooth family of symplectomorphisms $\upphi'_t:V\to V$ (for $t \in [0,1])$, where $\upphi'_0=id_V$ and $\upphi'_1 = \upphi \circ \upphi'$.
\end{defn}

Note that: $U,V \in C_{\X}$, each equipped with a symplectic form as part of their $\G$-structure, and $\upphi$ is a symplectomorphism serving as a weak equivalence in the context of symplectic homotopy.

\begin{defn}
    A map $\upphi : U \to V$ in $C_{\X}$ is a symplectic equivalence if it is a symplectomorphism, and it belongs to a class $W$.(e.g. isotopic to the identity or inducing isomorphisms on symplectic invariants.)
\end{defn}

\begin{defn}
    A sheaf $F$ in the category $\shv(C_{\X})$ is a symplectic motivic space if, for every $\upphi:U\to V$ in $W$, the map $F(\upphi): F(V) \to F(U)$ is an equivalence.
\end{defn}
\newline
\newline 
Form the category of spectrum objects $\Sp(\shv(C_{\X}))$ by stabilizing $\shv(C_{\X})$. This involves introducing spectrum objects, which are sequences of objects in $\shv(C_{\X})$ with bonding maps, similar to spectra in stable homotopy theory.

\begin{defn} 
From now on, we define $\mot_{\symp}(\X)$ as a subcategory of $\Sp(\shv(C_{\X}))$. An object is a pair $(E,\omega_E)$, where $E \in \Sp(\shv(C_{\X}))$ is a spectrum object, and $\omega_E: E \otimes E  \to \theta_{\X}[d]$ is a non-degenerate, skew-symmetric symplectic form, $\theta_{\X}$ is the global Thom spectrum object in $\Sp(\X)$ and $d$ a degree parameter. Maps in $\mot_{\symp}(\X)$ are morphisms in $\Sp(\shv(C_{\X}))$ that preserve the symplectic form up to homotopy.
\end{defn}
Since $\Sp(\shv(C_{\X}))$ is stable, ensure $\mot_{\symp}(\X)$ is closed under suspensions, loops, and cofiber sequences, while preserving the symplectic structure.

\begin{defn}
    The stabilized category of symplectic motives, denoted \(\mot_{\symp}^{\mathrm{stab}}(\X)\), is defined as follows:
    \begin{enumerate}[label=(\arabic*)]
         \item {Objects}: Triples \((E, q, \omega_E)\), where: \(E \in \Sp(\shv(C_{\X}))\) is a spectrum object, \(q \in \mathbb{Z}\) is the symplectic weight, \(\omega_E : E \otimes E \to \theta_{\X}[2q]\) is a non-degenerate, skew-symmetric symplectic form.
         
        \item {Morphisms}: A morphism from \((E, q_E, \omega_E)\) to \((F, q_F, \omega_F)\) is a map \(f : \Sigma^{q_F - q_E} E \to F\) in \(\Sp(\shv(C_{\X}))\) such that the following diagram commutes up to homotopy:
        \begin{equation*}
            \begin{CD}
                \Sigma^{q_F - q_E} E \otimes \Sigma^{q_F - q_E} E @>{f \otimes f}>> F \otimes F \\
                @V{\Sigma^{2(q_F - q_E)} \omega_E}VV @VV{\omega_F}V \\
                \theta_{\X}[2 q_F] @>{\cong}>> \theta_{\X}[2 q_F]
            \end{CD}
        \end{equation*}
        where the bottom map is the canonical isomorphism induced by the suspension.
        
        \item {Stable Structure}:
        \begin{itemize}
            \item {Suspension Functor}: \(\Sigma : \mot_{\symp}^{\mathrm{stab}}(\X) \to \mot_{\symp}^{\mathrm{stab}}(\X)\) defined by:
            \begin{equation*}
                \Sigma(E, q, \omega_E) = (\Sigma E, q + 1, \omega_{\Sigma E})
            \end{equation*}
            where \(\omega_{\Sigma E} : \Sigma E \otimes \Sigma E \theta_{\X}[2(q+1)]\) is explicitly given by the composition:
            \begin{equation*}
                \omega_{\Sigma E} : \Sigma E \otimes \Sigma E \cong \Sigma^2 (E \otimes E) \xrightarrow{\Sigma^2 \omega_E} \Sigma^2 \theta_{\X}[2q] \cong \theta_{\X}[2(q + 1)].
            \end{equation*}
            
            \item {Loop Functor}: \(\Omega(E, q, \omega_E) = (\Omega E, q - 1, \omega_{\Omega E})\), with \(\omega_{\Omega E}\) defined analogously using the loop functor.
            
            \item {Exact Triangles}: Inherited from \(\Sp(\shv(C_{\X}))\), with the symplectic structure preserved through the morphisms.
        \end{itemize}
    \end{enumerate}
\end{defn}
It contains $\mot_{\symp}(\X)$ for fixed $q$. 

\subsection{Symplectic motivic cohomology}
\begin{defn}
    For each $r \in \mathbb{Z}$, define symplectic motivic sphere: 
    \begin{equation*}
        \mathbb{S}_{\symp}^r =(\theta_{\X}[0],r,\omega_r), 
    \end{equation*}
    where $\theta_{\X}$ is the Thom spectrum at degree 0, $\omega_r:\theta_{\X}[0] \otimes \theta_{\X}[0] \to \theta_{\X}[2r]$ is defined as the composition:
    \begin{equation*}
        \omega_r : \theta_{\X} \otimes \theta_{\X} \xrightarrow[]{\omega} \theta_{\X}[2]  \xrightarrow[]{\Sigma^{2(r-1)}} \theta_{\X}[2r],
    \end{equation*}
    with $\omega$ being a canonical symplectic form on $\theta_{\X}$, and $\Sigma^{2(r-1)}$ the $2(r-1)$-fold suspension(or loop map if $r < 1$).
\end{defn}

\begin{defn}
    Define functor $\Sigma_{\omega}: \mot_{\symp}^{\mathrm{stab}}(\X) \to \mot_{\symp}^{\mathrm{stab}}(\X)$ that shifts the weight 
    \begin{equation*}
        \Sigma_{\omega}(E,q,\omega)=(E,q+1,\omega'),
    \end{equation*}
    where the new symplectic form  $\omega'$ is:
    \begin{equation*}
        \omega' = \omega_E[2]:E \otimes E \xrightarrow[]{\omega_E} \theta_{\X}[2q] \xrightarrow[]{\sum^2} \theta_{\X}[2(q+1)].
    \end{equation*}
\end{defn}
    
\begin{defn}
For integers \( p, q \in \mathbb{Z} \), define the bigraded suspension functor as follows:
\[
\Sigma^{p, q}\left(E, q_E, \omega_E\right) = \left\langle \Sigma^{p} E, q_E + q, \omega_{p, q} \right\rangle,
\]
where:
\begin{itemize}
    \item \( \Sigma^{p} E \) denotes the \( p \)-fold suspension of \( E \) in the stable \(\infty\)-category \(\operatorname{Sp}(\operatorname{Shv}(C_{\X}))\).
    \item \( q_E + q \) represents the new symplectic weight.
    \item \( \omega_{p, q}: \Sigma^{p} E \otimes \Sigma^{p} E \rightarrow \theta_{\X}[2(q_E + q)] \) is the adjusted symplectic form, given by \(\omega_{p, q} = \Sigma^{p} \omega_E \circ \tau,\)
    with:
    \begin{itemize}
        \item \( \tau: \Sigma^{p} E \otimes \Sigma^{p} E \rightarrow \Sigma^{p} (E \otimes E) \), the canonical isomorphism in the symmetric monoidal category \(\operatorname{Sp}(\operatorname{Shv}(C_{\X}))\).
        \item \( \Sigma^{p} \omega_E: \Sigma^{p} (E \otimes E) \rightarrow \Sigma^{p} \theta_{\X}[2 q_E] \), where we use the identification \( \Sigma^{p} \theta_{\X}[2 q_E] \cong \theta_{\X}[2 q_E + p] \).
        \item The condition \( p = 2 q \) ensures degree alignment, so that \( \theta_{\X}[2 q_E + p] = \theta_{\X}[2(q_E + q)] \).
    \end{itemize}
\end{itemize}
\textbf Note that: The functor \(\Sigma^{p, q}\) is defined only when \( p = 2 q \), ensuring that the symplectic form’s target matches the required degree \( \theta_{\X}[2(q_E + q)] \).
\end{defn}

\begin{defn}
    Construct a bigraded category $\mot_{\symp}^{\mathrm{bistab}}(\X)$:
    \begin{enumerate} [label=(\arabic*)]
        \item Objects: Same as $\mot_{\symp}^{\mathrm{stab}}(\X)$, triples (E,q,$\omega$).
        \item Morphisms: For ($E,q_E,\omega_E)$ and ($F,q_F,\omega_F)$, the mapping space is: 
            \begin{equation*}
                \map_{\mathrm{bistab}}((E,q_E,\omega_E),(F,q_F,\omega_F)) = \map_{\Sp(\shv(C_{\X}))}(\Sigma^{q_F-q_E}E,F),
            \end{equation*}
            with maps $f: \Sigma^{q_F-q_E}E \to F$ such that the symplectic forms are compatible, i.e. $\omega_F \circ (f \otimes f) \cong \Sigma^{2(q_F-q_E)}\omega_E$ under the suspension adjustment. Composition defined via the stable structure of $\Sp(\shv(C_{\X}))$.
    \end{enumerate}
\end{defn}
This category allows weight shifts via suspension, making it suitable for bigraded cohomology.

\begin{defn}
    The mapping space is: 
    \begin{equation*}
        [\Sigma^{p,q}(E,q_E,\omega_E),\mathbb{S}_{\symp}^0] = \map_{\Sp(\shv(C_{\X}))}(\Sigma^{p-q_E-q}E,\theta_{\X}[0])
    \end{equation*}
    with the symplectic condition enforced.
\end{defn}

\begin{defn}
    The symplectic motivic cohomology groups are defined as the connected components of this mapping space:
    \begin{equation*}
        H_{\mathrm{\symp}}^{p,q}((E,q_E,\omega_E),\mathbb{Z}) = \pi_0[ \Sigma^{p,q}(E,q_E,\omega_E),\mathbb{S}_{\symp}^0 ],
    \end{equation*}
    Substituting the mapping space:
    \begin{equation*}
     H_{\mathrm{\symp}}^{p,q}((E,q_E,\omega_E),\mathbb{Z}) = \pi_0 \map_{\Sp(\shv(C_{\X}))}(\Sigma^{p-q_E-q}E,\theta_{\X}[0]),
    \end{equation*}
    here $\pi_0$ takes the set of path components, giving a group (with $\mathbb{Z}$ indicating integer coefficients implicitly via the Thom spectrum).
\end{defn}
\\
What I'm going to say now will be covered in detail in the next paper. Assuming $\X = Shv(\mathcal{Y})$.
\\
Consider $u: C_{\X} \to Nis(\mathcal{Y})$ mapping $(U,\mathcal{O}_{U}) \to U$, where $U$ is regarded as an object in $Nis(\mathcal{Y})$ via its étale morphism to $\mathcal{Y}$. And $u(\upphi) = \upphi : V \to W$, a morphism in $Nis(\mathcal{Y})$. Since $C_{\X}$ has more coverings, the pushforward: 
 \begin{equation*}
u_* : PSh(C_{\X}) \to PSh_{Nis}(\mathcal{Y}), u_*F(U) = F((U,\mathcal{O}_U)),
  \end{equation*}
To ensure the result is a Nisnevich sheaf, compose with the Nisnevich sheafification functor $a_{Nis}$:
 \begin{equation*}
u_*^{sh} : Shv(C_{\X}) \to Shv_{Nis}(\mathcal{Y}), u_*^{sh}F = a_{Nis} \circ u_*F,
  \end{equation*}
This $u_*^{sh}$ maps sheaves on $C_{\X}$ to Nisnevich sheaves on $\mathcal{Y}$, resolving the topology mismatch. This extends to spectra:
\begin{equation*}
Sp(u_*^{sh}) : Sp(Shv(C_{\X})) \to Sp(Shv_{Nis}(\mathcal{Y}))
\end{equation*}
mapping $E  \mapsto  u_*^{sh}E.$ The forgetful functor $F: Mot_{symp}^{bistab}(\X) \to DM(\mathcal{Y})$ (Voevodsky's category of motives over $\mathcal{Y}$) is $F((E,q,\omega_E)) = L_{\mathbb{A}^1} \circ u_*^{sh}E$, where $L_{\mathbb{A}^1}$ is $\mathbb{A}^1$- localization and mapping morphism $f$ (respecting symplectic structure up to homotopy) directly. This functor forgets the symplectic structure and weight, mapping the spectrum $E$ to motive in $DM(\mathcal{Y})$ via the relationship between $C_{\X}$ and $Nis(\mathcal{Y})$, assuming $\X = Shv(\mathcal{Y})$ and $\mathcal{Y}$ is the underlying scheme.
\\
\\
Note that: This paper is just memos. We’ll submit the completed paper next time. I don't know when is it.(I don't have time because of my job.) If you have any questions, please feel free to contact me: Eita\_Haibara@protonmail.com
\\
\\
Thank you so much Taewan Kim for providing motivation and helping me with latex work(I'm not good at latex, so it was a great help). I would also like to thank Shelly Miyano for supporting the $\G$-scheme research. 

\bibliographystyle{unsrt}

\begin{thebibliography}{1}

\bibitem{haibara1}
Eita Haibara and Shelly Miyano.
\newblock Some memos on stable symplectic structured space, 2025.

\bibitem{grothendeick57}
A.~Grothendieck.
\newblock Sur la classification des fibres holomorphes sur la sphere de riemann.
\newblock {\em American Journal of Mathematics}, 79(1):121--138, 1957.

\bibitem{Blumberg_2013}
Andrew~J Blumberg, David Gepner, and Gonçalo Tabuada.
\newblock A universal characterization of higher algebraic k-theory.
\newblock {\em Geometry \& Topology}, 17(2):733–838, April 2013.

\bibitem{ha}
Jacob Lurie.
\newblock Higher algebra.
\newblock Unpublished. Available online at \url{https://www.math.ias.edu/~lurie/}, 09 2017.

\end{thebibliography}

\Large{Visiting our page! : \href {https://sites.google.com/view/pocariteikoku/home}{https://sites.google.com/view/pocariteikoku/home}}

\end{document}